# An algebraic approach for modeling and simulation of road traffic networks

Nadir Farhi[*,1], Habib Haj-Salem[1]
and Jean-Patrick Lebacque[1]

[1] *Université Paris-Est, IFSTTAR/COSYS/GRETTIA, F-77447 Champs-sur Marne Cedex France*

**Abstract.** We present in this article an algebraic approach to model and simulate road traffic networks. By defining a set of road traffic systems and adequate concatenating operators in that set, we show that large regular road networks can be easily modeled and simulated. We define elementary road traffic systems which we then connect to each other and obtain larger systems. For the traffic modeling, we base on the LWR first order traffic model with piecewise-linear fundamental traffic diagrams. This choice permits to represent any traffic system with a number of matrices in specific algebraic structures. For the traffic control on intersections, we consider two cases: intersections controlled with a priority rule, and intersections controlled with traffic lights. Finally, we simulate the traffic on closed regular networks, and derive the macroscopic fundamental traffic diagram under the two cases of intersection control.

**Keywords:** Road traffic modeling and simulation, min-plus algebra, traffic control.

## 1 Introduction

Modeling the traffic in urban networks is necessary to understand the vehicular dynamics and set adequate strategies and controls, in order to improve the service. Many models with different approaches exist in the literature (1). We present in this article a urban traffic model based in the cell-transmission model (2) (a numerical scheme of the first order macroscopic LWR model (3), (4)); see also (5). The model adapts the existing approach to the urban traffic framework. Moreover, two models of intersection control are proposed. An algebraic formulation of the whole vehicular dynamics in a urban road network is made. The formulation permits to represent the traffic dynamics in the network by a number of matrices in the min-plus algebra (a specific algebraic structure) (6).

The approach we adopt here is a system theory approach, where the urban traffic network is build from predefined elementary traffic systems and adequate operators, for the connection of these systems. We first present the link traffic model inspired from the cell-transmission model (2), with its algebraic formulation. In section 3, we give two intersection traffic models. In section 4, we explain the algebraic construction of an American-like (regular) city, by giving the three elementary traffic systems and the main operator we use for that. In the last part of the article we present some numerical traffic simulations on regular cities set on a torus (closed networks). This configuration permits to easily derive the macroscopic fundamental diagram on such networks. Finally, we discuss the traffic phases obtained from those diagrams, under two control policies set on the intersections. In this article, we only review the traffic models we use. For more details on those models; see (8), (9), (10) and (12). The main contribution of this article is the system theory approach we propose for building and simulating urban traffic networks.

---

[*] Corresponding author (`nadir.farhi@ifsttar.fr`)

## 2  The link model

The model we propose here is based on the macroscopic first order LWR model (3) (4), with triangular fundamental diagrams, where the dynamics of vehicle pelotons moving through road sections is described. We assume here that only pelotons are observed. Moreover, the density of pelotons is considered to be binary, in the sense that, at a given time, the density on a given section is equal to 1 if any peloton of vehicles is moving on, and it is equal to zero otherwise. We think that this mesoscopic representation of the traffic dynamics is convenient to describe the traffic in urban networks.

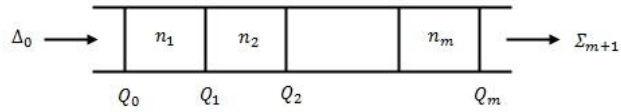
Figure 1. A single-lane road.

We first present the traffic model on a single link, where traffic is unidirectional, and where vehicles move without passing. Let us explain the traffic dynamics on a road of m sections. We use the following notations.
- $n_i(t) \in \{0,1\}$: number of pelotons in section $i$, at time $t$, with $i = 1,2,\ldots,n$ and $t \in \mathbb{N}$.
- $\bar{n}_i(t) = 1 - n_i(t) \in \{0,1\}$: free space in section $i$ at time $t$, with $i = 1,2,\ldots,n$ and $t \in \mathbb{N}$.
- $Q_0(t) \in \mathbb{N}$: Cumulated flow (in number of pelotons per time unit) from time zero to time $t$, of vehicle pelotons entering to section 1.
- $Q_i(t) \in \mathbb{N}, i = 1,2,\ldots,n-1$: Cumulated flow from time zero to time $t$, of vehicle pelotons passing from section $i$ to section $i+1$.
- $Q_m(t) \in \mathbb{N}$: Cumulated flow from time zero to time $t$, of vehicle pelotons leaving section $m$.

We assume here triangular fundamental diagrams on all the road sections.

$$q_i = \min\{v_i\, \rho_i, w_i\left(\rho_i^j - \rho_i\right)\}. \tag{1}$$

where $q_i, \rho_i, v_i, w_i$ and $\rho_i^j$ denote respectively, the car-flow, the car-density, the free speed, the backward wave speed, and the jam density, in section $i$. We assume that all sections have the same fundamental diagram. Moreover, according to the assumptions above, we assume that $v_i = w_i = \rho_i^j = 1, \forall i = 1,2, \ldots, m$. We thus obtain the following fundamental diagram for all the sections.

$$q_i = \min\{\rho_i, (1 - \rho_i)\} \tag{2}$$

According to the cell-transmission model (2) (7), which is a convenient numerical scheme of the LWR macroscopic model (3) (4), the traffic demand and supply are derived from the fundamental traffic diagram, and are given as follows.

- $\delta_i(t) = \min(v_i \rho_i(t), q_i^{max}) = \min(\rho_i(t), 1/2)$: the traffic demand from section $i$ to section $i+1$ at time t.
- $\sigma_i(t) = \min(q_i^{max}, w_i\left(\rho_i^j - \rho_i(t)\right)) = \min(1/2, 1 - \rho_i(t))$: the traffic supply of section $i$ to section $i-1$.

where

$$q_i^{max} = \frac{\rho_j}{\frac{1}{v_i} + \frac{1}{w_i}} = \frac{1}{2}, \forall i = 1,2, \ldots, m. \tag{3}$$

The *cumulated* traffic demand in the entry of the road, denoted by $\Delta_0(t)$, as well as the *cumulated* traffic supply on the exit of the road, denoted by $\Sigma_{n+1}(t)$ are supposed to be given over the whole time. They represent the boundary conditions of the system. The initial traffic condition consists here in giving the densities $\rho_i(0), i = 1,2, \ldots, m$ (the densities on each road section at time zero).

We assume that all the sections of the road have the same length, which we denote by $\Delta x$. Moreover, we fix the time unit $dt$ to $dt = \Delta x/v = \Delta x/w$. The model consists finally in giving the dynamics of the cumulated flows $Q_i(t), i = 0,1, \ldots, m$ over time $t \in \mathbb{N}$.

$$\begin{aligned} Q_0(t+dt) &= \min\{\Delta_0(t), Q_1(t) + \bar{n}_1(0)\} \\ Q_i(t+dt) &= \min\{Q_{i-1}(t) + n_i(0), Q_{i+1}(t) + \bar{n}_{i+1}(0)\} \\ Q_m(t+dt) &= \min\{Q_{m-1}(t) + n_m(0), \Sigma_{n+1}(t)\} \end{aligned} \tag{4}$$

and, by that, updating the number of pelotons $n_i(t), i = 1,2, \ldots, n; t \in \mathbb{N}$.

$$n_i(t) = n_i(0) + Q_{i-1}(t) - Q_i(t), \qquad i = 1,2, \ldots, n. \tag{5}$$

Let us notice that we assume here that the cumulated flows are initialized to zero: $Q_i(0) = 0, \forall i = 0,1,\ldots,m$.

**Algebraic formulation**

We consider here the algebraic structure $\mathbb{R}_{min} := (\mathbb{R} \cup \{+\infty\}, \oplus, \otimes)$, where the operations $\oplus$ and $\otimes$ are defined as follows.

$$a \oplus b := \min(a,b), \quad \forall a,b \in \mathbb{R}_{min}$$
$$a \otimes b := a + b, \quad \forall a,b \in \mathbb{R}_{min}$$

The structure $\mathbb{R}_{min}$ is a dioid (an idempotent semiring); see (6). We denote by $\varepsilon = +\infty$ and $e = 0$ respectively the zero and the unity elements for $\mathbb{R}_{min}$. We have also a dioid in the set $\mathcal{M}_{n \times n}(\mathbb{R}_{min})$ of square matrices with elements in $\mathbb{R}_{min}$, where the operations $\oplus$ and $\otimes$ are defined as follows.

$$(A \oplus B)_{ij} := A_{ij} \oplus B_{ij} = \min(A_{ij}, B_{ij}), \quad \forall A,B \in \mathcal{M}_{n \times n}(\mathbb{R}_{min})$$
$$(A \otimes B)_{ij} := \bigoplus_{1 \leq k \leq n} (A_{ik} \otimes B_{kj}) = \min_{1 \leq k \leq n}(A_{ik} + B_{kj}), \quad \forall A,B \in \mathcal{M}_{n \times n}(\mathbb{R}_{min}).$$

It is then easy to check that the dynamics (4) can be written as follows.

$$Q(t + dt) = A \otimes Q(t) \oplus b(t) \qquad (6)$$

where $Q(t)$ is the vector whose components are the cumulated flows $Q_i(t)$, and where $A \in \mathcal{M}_{n \times n}(\mathbb{R}_{min})$ and $b(t) \in \mathcal{M}_{1 \times n}(\mathbb{R}_{min})$ are given as follows.

$$A = \begin{pmatrix} \varepsilon & \bar{n}_1(0) & \varepsilon & \cdots & & \cdots & \varepsilon \\ n_1(0) & \varepsilon & \bar{n}_2(0) & \varepsilon & & \cdots & \varepsilon \\ \varepsilon & n_2(0) & \varepsilon & \bar{n}_3(0) & & & \varepsilon \\ \vdots & & \ddots & \ddots & \ddots & \ddots & \vdots \\ & & & n_{m-1}(0) & \varepsilon & \bar{n}_m(0) \\ & & & & n_m(0) & \varepsilon \end{pmatrix}, \quad b(t) = \begin{pmatrix} \Delta_0(t) \\ \varepsilon \\ \varepsilon \\ \vdots \\ \varepsilon \\ \Sigma_{n+1}(t) \end{pmatrix}.$$

with $Q(0) = 0$.

With this formulation, the traffic model on any single-lane road is summarized by the two matrices $A$ and $b(t)$, $t \in \mathbb{N}$. The simulation of the traffic model is then simply done by iterating the min-plus linear dynamics (6), with the initial condition $Q(0) = 0$. We notice that the matrix $A$ and the vector $b(t)$ contain respectively the initial condition (initial density) and the boundary conditions (demand inflow and supply outflow). For more details on the model presented in this section, see (8) (9). We will see below (in the two dimensional traffic modeling section), that the linearity of the traffic dynamics obtained in the one dimension model cannot be preserved.

## 3 Two dimensional traffic modeling

In order to be able to model the traffic on road networks, we need to have models for intersections. We present in this section two models. The first model describes the traffic inflowing to and out-flowing from an intersection with two entry roads and two exit roads where one of the entry roads has priority with respect to the other one. The second model considers that the intersection is controlled with a traffic light.

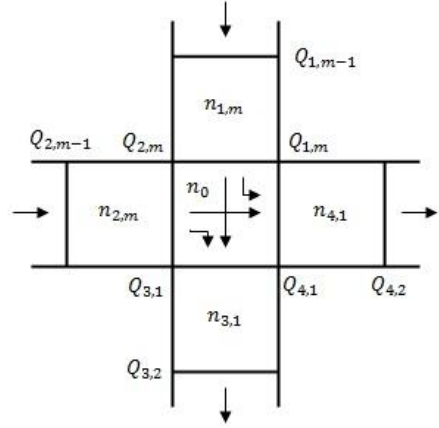

Figure 2. Intersection of two roads.

### 3.1 Intersection model with a priority rule.

Let us consider the intersection of Figure 2, where a priority rule is set. Vehicles entering the intersection from road 1 (the North) have priority with respect to vehicles entering the intersection from road 2 (the West). $n_0(t)$ and $\bar{n}_0(t)$ denote respectively the number of pelotons and the free space in the intersection at time $t$. Equations (7) below only describe the traffic dynamics on the intersection. The traffic on the roads follows the dynamics described above.

$$\begin{aligned}
Q_{1m}(t + dt) &= \min\{Q_{1m-1}(t) + n_{1m}(0), Q_{31}(t) + Q_{41}(t) - Q_{2m}(t) + \bar{n}_0(0)\} \\
Q_{2m}(t + dt) &= \min\{Q_{2,m-1}(t) + n_{2,m}(0), Q_{31}(t) + Q_{41}(t) - Q_{1m}(t + dt) + \bar{n}_0(0)\} \\
Q_{31}(t + dt) &= \min\{\alpha_{13}Q_{1m}(t) + \alpha_{23}Q_{2m}(t) + n_0(0), Q_{32}(t) + \bar{n}_{31}(0)\} \\
Q_{41}(t + dt) &= \min\{\alpha_{14}Q_{1m}(t) + \alpha_{24}Q_{2m}(t) + n_0(0), Q_{42}(t) + \bar{n}_{41}(0)\}
\end{aligned} \quad (7)$$

where the notations used in (7) are (see Figure 2):
- $Q_{1m}(t)$: cumulated outflow from road 1, which is also the cumulated inflow to the intersection, from the north side, up to time $t$.
- $Q_{2m}(t)$: cumulated outflow from road 2, which is also the cumulated inflow to the intersection, from the west side, up to time $t$.
- $Q_{31}(t)$: cumulated outflow from the intersection to the south, which is also the cumulated inflow to road 3, up to time $t$.

- $Q_{41}(t)$: cumulated outflow from the intersection to the east, which is also the cumulated inflow to road 4, up to time $t$.

The dynamics of $Q_{1m}$ and $Q_{2m}$ in (7) (the two first equations) set the priority to the outflow from road 1 with respect to the outflow from road 2. This is done by the introduction of an implicit term in the dynamics of $Q_{2m}$ in (7). For more details, see (8) (9) (10).

Using the same notations as above, we can easily check that the dynamics (7) is written with matrix notations as follows.

$$Q(t + dt) = D \otimes \left( H\, Q(t) + G\, Q(t + dt) \right) \oplus b(t) \qquad (8)$$

where $D$ is a min-plus matrix, and $H$ and $G$ are standard matrices. The matrices $H$ and $G$ contain multipliers that cannot be expressed linearly in the min-plus algebra. These multipliers are needed to model the turning rates as well as the priority rule at the intersection. The turning rates in the level of the intersection are given by $H$ and $G$, where $H$ gives the turning rates with a time delay $dt$, and $G$ gives the turning rates without any time delay (the implicit term setting the priority rule). For more details on the model presented in this section, see (8) (9) (10).

### 3.2 Intersection model with a traffic light control.

We give in this section the traffic dynamics in the case where the intersection is managed by means of a traffic light. In a first step, we consider only the case where an open loop control is set on the traffic light. The control is assumed to be periodic with a time period (cycle) denoted by $c$ (which is in fact equal to $c\, dt$). The green times for the north and the west sides are denoted respectively by $g_N$ and $g_W$. The integral red times between the two green times are denoted by $r_1$ and $r_2$ respectively for the integral red time from the end of $g_N$ and the beginning of $g_W$ and for the integral red time from the end of $g_W$ and the beginning of $g_N$; see Figure 3.

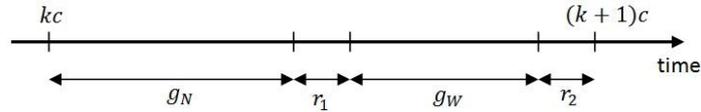

Figure 3. Time cycle for the traffic light.

The dynamics (7) is modified to:

$$Q_{1m}(t+dt) = \min \begin{cases} Q_{1m-1}(t) + n_{1m}(0), \\ Q_{31}(t) + Q_{41}(t) - Q_{2m}(t) + \bar{n}_0(0), \\ Q_{1m}(t) + L_1. \end{cases}$$

$$Q_{2m}(t+dt) = \min \begin{cases} Q_{2,m-1}(t) + n_{2,m}(0), \\ Q_{31}(t) + Q_{41}(t) - Q_{1m}(t) + \bar{n}_0(0), \\ Q_{2m}(t) + L_2. \end{cases} \quad (9)$$

$$Q_{31}(t+dt) = \min\{\alpha_{13}Q_{1m}(t) + \alpha_{23}Q_{2m}(t) + n_0(0), Q_{32}(t) + \bar{n}_{31}(0)\}$$
$$Q_{41}(t+dt) = \min\{\alpha_{14}Q_{1m}(t) + \alpha_{24}Q_{2m}(t) + n_0(0), Q_{42}(t) + \bar{n}_{41}(0)\}$$

where $L_1 = \begin{cases} q_{1,m}^{max} = \frac{1}{2} & \text{if } t \in [kc, kc + g_N], \\ 0 & \text{otherwise} \end{cases}$

and $L_2 = \begin{cases} q_{2,m}^{max} = \frac{1}{2} & \text{if } t \in [kc + g_N + r_1, kc + g_N + r_1 + g_W], \\ 0 & \text{otherwise} \end{cases}$

Thus, in the time instants when $L_1 = q_{1,m}^{max} = 1/2$, the traffic light is green for the road 1, because, $Q_{1,m}(t)$ may be increased by $q_{1,m}^{max}$, under the two constraints of upstream demand and downstream supply. In the time instants when $L_1 = 0$, the traffic light is red for road 1, because, $Q_{1,m}(t)$ stays constant, i.e. $Q_{1,m}(t+dt) = Q_{1,m}(t)$. The same reasoning is made for the road 2. The algebraic formulation of the model (9) is similar to the one done in (8), but we need here to define four dynamics, one for each phase of the time cycle. For more details in the model presented in this section, see (8) (10) (9).

## 4  An American-like city

We define in this section a set of dynamic systems such that any traffic system defined under the models presented above, is contained in that set. We also define operators for the connection of those systems. The systems we consider here are those with two vectors of input signals $U$ and $V$, two vectors of state signals $P$ and $Q$, and two vectors of output signals $Y$ and $Z$, such that we can write

$$\begin{pmatrix} P(t+dt) \\ Q(t+dt) \\ Y(t+dt) \\ Z(t+dt) \end{pmatrix} = \begin{pmatrix} 0 & A & 0 & B \\ C & \varepsilon & D & \varepsilon \\ 0 & E & 0 & 0 \\ F & \varepsilon & \varepsilon & \varepsilon \end{pmatrix} \boxtimes \begin{pmatrix} P(t+dt) \\ Q(t) \\ U(t+dt) \\ V(t) \end{pmatrix}$$

$$:= \begin{pmatrix} AQ(t) + BV(t) \\ C \otimes P(t+dt) \oplus D \otimes U(t+dt) \\ EQ(t) \\ F \otimes V(t) \end{pmatrix}, \quad (10)$$

where $A, B$ and $E$ are standard matrices, while $C, D$ and $F$ are min-plus matrices. This construction is inspired from Petri Net modeling, see (8). If we denote by $S$ the system (10), then we write $(Y, Z) = S(U, V)$. Let us explain how traffic dynamics given above are written in the form (10). For that, we first do it for the three elementary systems on which we will base for building traffic systems of large networks. The three elementary systems that we consider here are the following.

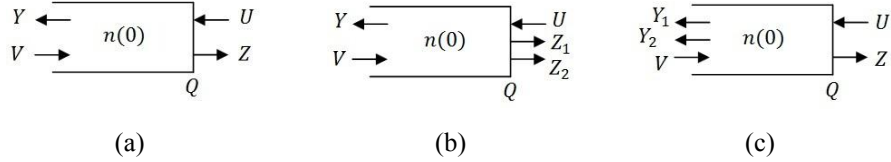

(a)  (b)  (c)

Figure 4. Elementary traffic systems: (a) a road section, (b) an intersection entry, (c) an intersection exit.

a) a road section is the elementary traffic system in a road. The system has two input signals $U$ and $V$, one state signal $Q$, and two output signals $Y$ and $Z$.
b) an intersection entry is a special road section with more output signals than an ordinary road section (a). The system has two input signals $U$ and $V$, one state signal $Q$, and three output signals $Y, Z_1$ and $Z_2$.
c) an intersection exit is a special road section with more output signals than an ordinary road section (a). The system has two input signals $U$ and $V$, one state signal $Q$, and three output signals $Y_1, Y_2$ and $Z$.

In order to clarify how the dynamics of these elementary systems are written in the form (10), we explain the dynamics of a road section (system (a)). Following the dynamics (4) (or (6)), the dynamics of the road section (a) is written as follows.

$$Q(t + dt) = \min(n(0) + V(t), U(t + dt)), \qquad (11)$$
$$Y(t + dt) = Q(t + dt),$$
$$Z(t + dt) = \bar{n}(0) + Q(t).$$

Then by introducing intermediate variables, we get

$$P_1(t + dt) = V(t + dt)$$
$$P_2(t + dt) = Q(t + dt)$$
$$Q(t + dt) = \min(n(0) + P_1(t), U(t + dt)), \qquad (12)$$
$$Y(t + dt) = Q(t + dt),$$
$$Z(t + dt) = \bar{n}(0) + P_2(t).$$

Which can be easily written in the form (10) with

$$P = \begin{pmatrix} P_1 \\ P_2 \end{pmatrix}, A = \begin{pmatrix} 0 \\ 1 \end{pmatrix}, B = \begin{pmatrix} 1 \\ 0 \end{pmatrix}, C = (n(0) \quad \varepsilon), D = e, E = 1, F = (\varepsilon \quad \bar{n}(0)).$$

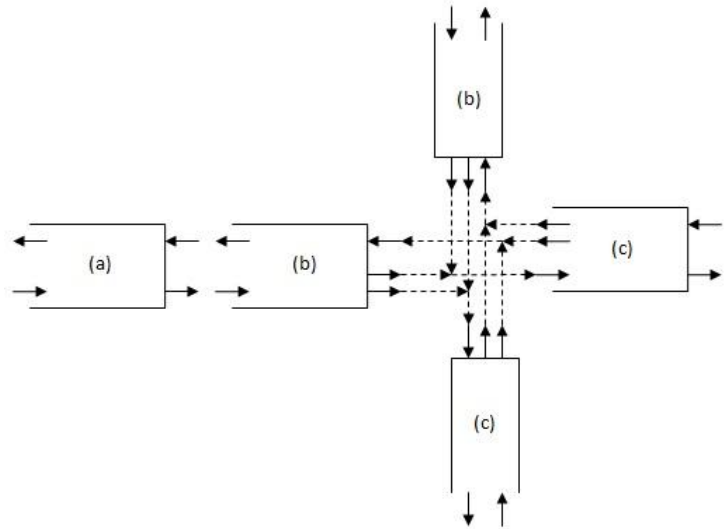

Figure 5. Connection of traffic elementary systems.

The dynamics of the two systems (b) and (c) are obtained in the similar way. Let us now explain how the systems are connected. For this, we define below the operator used for the connection. In figure 5, we illustrate the connection of road sections, and the construction of an intersection. Let us notice that an intersection is composed of two intersection entries and two intersection exits.

**Connection of systems**

Connecting two system $S_1$ and $S_2$ consists in equaling a part of inputs of each system with a part of outputs of the other system. We thus need first to specify the parts of inputs and outputs to be equalized. Let us note $S_1{}_{Y,Z,Y',Z'}^{U,V,U',V'}$ and $S_2{}_{Y,Z'',U',V'}^{U,V,Y',Z'}$, where $U', V'$ are inputs for $S_1$, and outputs for $S_2$, while $Y', Z'$ are inputs for $S_2$ and outputs for $S_1$. The connection of the two systems $S_1$ and $S_2$, denoted simply by $S_1 S_2$, is the system $\hat{S}_{Y,Z,Y'',Z''}^{U,V,U'',V''}$ given as the solution, on $Y, Y'', Z, Z''$, of the system

$$\begin{cases}(YY', ZZ') = S_1(UU', VV') \\ (U'Y, V'Z) = S_2(Y'U'', Z'V'')\end{cases}$$

Then, if we partition the input matrices of both systems $S_1$ and $S_2$ as follows

$$[B_1 B'_1], \quad [B'_2 B''_2], \quad [D_1 D'_1], \quad [D'_2 D''_2]$$

and the output matrices of the systems as follows

$$\begin{bmatrix}E_1 \\ E'_1\end{bmatrix}, \quad \begin{bmatrix}E'_2 \\ E''_2\end{bmatrix}, \quad \begin{bmatrix}F_1 \\ F'_1\end{bmatrix}, \quad \begin{bmatrix}F'_2 \\ F''_2\end{bmatrix},$$

then the system $\hat{S}$ is given by the matrices $\hat{A}, \hat{B}, \hat{C}, \hat{D}, \hat{E}$ and $\hat{F}$

$$\hat{A} = \begin{pmatrix} A_1 & 0 & 0 & B'_1 \\ 0 & A_2 & B'_2 & 0 \\ E_1 & 0 & 0 & 0 \\ 0 & E'_2 & 0 & 0 \end{pmatrix}, \hat{B} = \begin{pmatrix} B_1 & 0 \\ 0 & B''_2 \\ 0 & 0 \\ 0 & 0 \end{pmatrix}, \hat{C} = \begin{pmatrix} C_1 & \varepsilon & \varepsilon & D'_1 \\ \varepsilon & C_2 & D'_2 & \varepsilon \\ F'_1 & \varepsilon & \varepsilon & \varepsilon \\ \varepsilon & F'_2 & \varepsilon & \varepsilon \end{pmatrix},$$

$$\hat{D} = \begin{pmatrix} D_1 & \varepsilon \\ \varepsilon & D''_2 \\ \varepsilon & \varepsilon \\ \varepsilon & \varepsilon \end{pmatrix}, \hat{E} = \begin{pmatrix} E_1 & 0 & 0 & 0 \\ 0 & E''_2 & 0 & 0 \end{pmatrix}, \hat{F} = \begin{pmatrix} F_1 & \varepsilon & \varepsilon & \varepsilon \\ \varepsilon & F''_2 & \varepsilon & \varepsilon \end{pmatrix}.$$

For more details on this construction see (8).

**Closed loop control.**

We present in this section the application of an existing centralized urban control strategy, which is called TUC (Traffic Urban Control), see (11). The objective here is to derive the macroscopic fundamental traffic diagram on a regular city, under this control strategy, and then compare it to the diagrams obtained under the open loop control presented above, and under the priority rule.

TUC strategy assumes given a nominal traffic state (vehicle densities on the roads and controls in intersections), and regulates the traffic in the urban network, around the nominal traffic state. Let us use the notations.

- $x_i(t)$: the number of vehicles moving on raod $i$ at time $t$.
- $\bar{x}_i$: nominal number of vehicles moving on road $i$.
- $u_i(t)$: outflow from road $i$ at time $t$.
- $\bar{u}_i$: nominal outflow from road $i$.

We then solve the following linear quadratic control problem.

$$\min_{u \in U} \sum_{t=0}^{+\infty} (x(t) - \bar{x})' Q(x(t) - \bar{x}) + (u(t) - \bar{u})' R(u(t) - \bar{u}) \quad \text{(13)}$$
$$(x(t + dt) - \bar{x}) = (x(t) - \bar{x}) + B(u(t) - \bar{u}).$$

For example, according to Figure 2, the dynamics of the number of vehicles moving on the road 4 is written

$$x_4(t + dt) = x_4(t) + \alpha_{14} u_1(t) + \alpha_{24} u_2(t) - u_4(t). \quad \text{(14)}$$

For more details on this approach, see (11) (8).

## 5  Simulation and derivation of macroscopic fundamental diagram

Following the models presented in the sections above, we build a regular city (an American-like city, where parallel horizontal avenues with alternated senses intersect parallel vertical avenues with alternated senses, see Figure 6). Without loss of generality, we assume here that the city is wholly symmetric, in the sense that all roads have the same length, the same fundamental traffic diagram, and all the turning rates are equal to 1/2. Because of the symmetry, the ideal control in this configuration would be to uniformly distribute the number of cars on the roads of the city.

The objective of considering closed networks, like a city on a torus (Figure 6) is to be able to fix the car density on the whole network, and then derive the asymptotic average car-flow on the whole network. Theoretical results on the existence and uniqueness of such asymptotic average flows, as well as their dependence on the initial average car-density in the network, can be found in (8) (12).

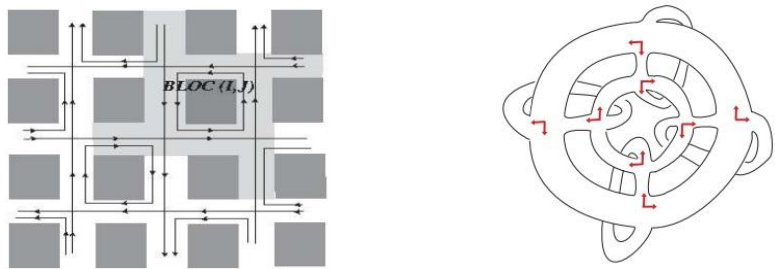

Figure 6. A regular city (left side), and a regular city on a torus.

In Figure 7 we give the fundamental traffic diagrams (average traffic flow in function of the average traffic density on the city) derived from the whole city on a torus, under different control strategies on the intersections.

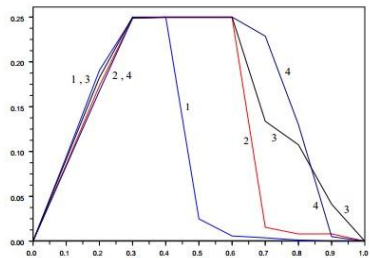

Figure 7. Comparison of the fundamental traffic diagrams obtained under different control policies set on the intersections of the regular city on a torus.
1. Priority rule. 2. Traffic lights in open loop, with equal green times for both directions of every intersection. 3. Local feedback that sets green times on every intersection, proportional to the densities on the two entering roads. 4. Centralized feedback control with TUC strategy.

In Figure 8, we show some simulations of traffic on the regular city on a torus. In particular, we compare in that figure the control of traffic lights under open loop and centralized closed loop controls. The result is that the centralized closed loop control is the better strategy, in the sense that it attains surely the nominal traffic state, which is here the uniform distribution of the number of cars on the roads of the city. This is also confirmed on the fundamental diagrams of Figure 7, where only the centralized feedback strategy guaranties acceptable flows in the case of high densities.

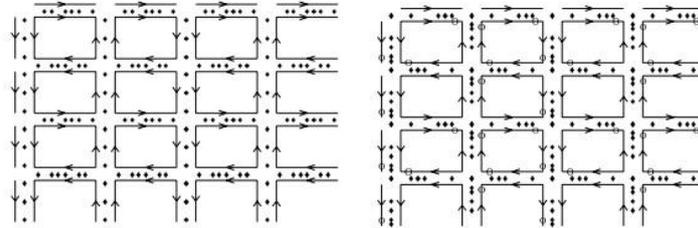

Figure 8. Traffic simulation. On the left side: open loop control. On the right side: centralized closed loop control.

## 6  Conclusion

The traffic modeling approach we proposed in this article permits to algebraically build large urban regular networks, such as American-like cities. Two intersection models are presented: intersection managed with a priority rule, and intersection controlled with a traffic light. Moreover, a centralized feedback control is applied to control such road networks. Finally we compared different control approaches by means of the derived macroscopic fundamental diagrams. The conclusion is that centralized feedback controls are the better control strategies for the stabilization of the traffic under severe congestion.